\documentclass[11 pt]{amsart}
\usepackage{amssymb}
\usepackage{amsmath}
\usepackage{amsfonts}
\usepackage{graphicx}
\usepackage{amsthm}
\usepackage{enumerate}
\usepackage[mathscr]{eucal}
\usepackage{verbatim}

 \theoremstyle{plain}
\newtheorem{theorem}{Theorem}[section]
\newtheorem{lemma}[theorem]{Lemma}

\numberwithin{equation}{section}
\theoremstyle{plain}

\theoremstyle{remark}

\def\bbR{{\mathbb {R}}}

\def\cS{\mathcal S}

\begin{document}

\date{May, 2007}

\title
{Restricted Radon transforms and projections of planar sets}

\author[]
{Daniel M. Oberlin}

\address
{D. M.  Oberlin \\
Department of Mathematics \\ Florida State University \\
 Tallahassee, FL 32306}
\email{oberlin@math.fsu.edu}

%\subjclass{42B10, 42B99}
%\keywords{Fourier transforms of measures on curves,
%Fourier restriction problem}

\thanks{
This work was supported in part by NSF grant DMS-0552041.
}

\begin{abstract}
We establish a mixed norm estimate for the Radon transform in $\bbR^2$ when the set of directions has fractional dimension. 
This estimate is used to prove a result about an exceptional set of directions connected with projections of planar sets. That leads to 
a conjecture analogous to a well-known conjecture of Furstenberg.
\end{abstract}

\maketitle

\section{Introduction}

For each $\omega\in S^1$, fix $\omega^\perp$ with $\omega^\perp \perp \omega$. Define a Radon transform $R$
for functions $f$ on $\bbR ^2$ by 
$$
Rf(t,\omega )=\int_{-1}^1 f(t\,\omega +s\,\omega^\perp )\,ds.
$$
Suppose $0<\alpha <1$ and fix a nonnegative Borel measure $\lambda$ on $S^1$ which is $\alpha$-dimensional in the sense 
that $\lambda (B(\omega ,\delta ))\lesssim \delta^\alpha$ for $\omega\in S^1$. 
We are interested in mixed norm estimates for $R$ of the following form: 
\begin{equation}\label{est}
\Big[\int_{S^1}\Big(\int_{-1}^1 |Rf(t,\omega )|^s dt\Big)^{q/s}d\lambda (\omega )\Big]^{1/q}\lesssim \|f\|_p .
\end{equation}
Here are some conditions which are necessary for \eqref{est}: testing on $f=\chi_{B(0,\delta )}$ shows that 
\begin{equation}\label{N1}
\frac{2}{p}\leq 1+\frac{1}{s};
\end{equation}
if there is $\omega_0 \in S^1$ such that $\lambda (B(\omega_0 ,\delta ))\gtrsim \delta^\alpha$ for small positive $\delta$, then
testing on $1$ by $\delta$ rectangles centered at the origin in the direction $\omega_0^\perp$  gives 
\begin{equation}\label{N2}
\frac{1}{p}\leq \frac{1}{s}+\frac{\alpha}{q};
\end{equation} 
if the Lebesgue measure in $S^1$ of the $\delta$-neighborhood in $S^1$ of the support of $\lambda$ is $\lesssim \delta^{1-\alpha}$, then 
testing on unions of $1$ by $\delta$ rectangles in the directions of the support of $\lambda$ gives 
\begin{equation}\label{N3}
\frac{1-\alpha}{p}\leq \frac{1}{s}.
\end{equation}
Our first result is that these necessary conditions are almost sufficient:
\begin{theorem}\label{RRTest}
Suppose $p,q,r\in [1,\infty ]$ satisfy the conditions \eqref{N1}, \eqref{N2}, and \eqref{N3} with strict inequality. Then the estimate 
\eqref{est} holds.
\end{theorem}

Now suppose that $\mu$ is a nonnegative Borel measure on $\bbR^2$. If $\omega\in S^1$, define the projection $\mu_{\omega}$ of $\mu$ in the direction of $\omega$ by 
$$ 
\int_{\bbR}f(y)\,d\mu_{\omega}(y)\doteq \int_{\bbR^2}f(x\cdot\omega )\,d\mu (x),
$$
where $x\cdot\omega$ denotes the inner product in $\bbR^2$. Fix $\alpha\in (0,1)$ and suppose that $\lambda$ is an $\alpha$-dimensional measure
on $S^1$. Then, for $\epsilon >0$, there is $C=C(\epsilon )$ such that
$$
\int_{S^1}\frac{d\lambda (\omega )}{|\omega \cdot \omega_0 |^{\alpha -\epsilon}}\leq C(\epsilon )
$$
for all $\omega_0 \in S^1$. 
The computation
$$
\int_{S^1}I_{\alpha -\epsilon}(\mu _{\omega})\,d\lambda (\omega)=
\int_{S^1}\int_{\bbR}\int_{\bbR}
\frac{d\mu_{\omega} (y_1 )d\mu_{\omega} (y_2 )}{|y_1 -y_2 |^{\alpha -\epsilon}}d\lambda (\omega )=
$$
$$
\int_{\bbR^2}\int_{\bbR^2}
\int_{S^1}
\frac{d\lambda (\omega )}
{|\omega \cdot 
\frac{x_1 -x_2 }{|x_1 -x_2 |}|^{\alpha -\epsilon}}
\frac{d\mu (x_1 )d\mu (x_2 )}{|x_1 -x_2 |^{\alpha -\epsilon}}\leq\,C(\epsilon )\,I_{\alpha -\epsilon}(\mu )
$$
is due to Kaufman \cite{K}. Refining an earlier result of Marstrand \cite{M}, it shows that if 
$E\subset\bbR^2$ has dimension $\beta \leq 1$ and $p_\omega (E)$ is the projection of $E$ onto the line through the origin in the direction of $\omega$, then 
\begin{equation}\label{dimest}
\dim\{\omega\in S^1 :\dim p_\omega (E) <\alpha\}\leq \alpha 
\end{equation}
whenever $\alpha\leq\beta$. (In this note \lq\lq$\dim$" stands for Hausdorff dimension.) In particular,
\begin{equation}\label{dimest2}
\dim\{\omega\in S^1 :\dim p_\omega (E) <\beta\}\leq \beta .
\end{equation}
The next theorem, whose analog for Minkowski dimension is trivial, complements Kaufman's results \eqref{dimest} and \eqref{dimest2}:
\begin{theorem}\label{result2}
If $\dim E =\beta \leq 1$ then 
\begin{equation}\label{dimest3}
\dim\{\omega\in S^1 :\dim p_\omega (E) <\beta /{2}\}=0 .
\end{equation}

\end{theorem}
\noindent
The estimates \eqref{dimest2} and \eqref{dimest3} lead naturally to the conjecture that if $\alpha\leq\beta\leq 1$ then
\begin{equation}\label{dimest4}
\dim\{\omega\in S^1 :\dim p_\omega (E) <(\alpha +\beta )/2\}\leq \alpha .
\end{equation}
One may view this conjecture as an analog of the conjecture that Furstenberg $\alpha$-sets have dimension
at least $(3\alpha +1)/2$, with \eqref{dimest} being the analog of the known $2\alpha$ lower bound 
for the dimension of Furstenberg sets  
and with \eqref{dimest3} being the analog of the known
$(\alpha +1)/2$ lower bound. Indeed, \eqref{dimest4} with $\beta =1$ would imply the Furstenberg conjecture for a certain class 
of model Furstenberg sets. (Information about Furstenberg's conjecture is contained in \cite{W}.)
The link between Theorems \ref{RRTest} and \ref{result2} is the fact that, formally,  
$\mu_{\omega}=R\mu (\cdot\, ,\omega )$.
\section{Proof of Theorem \ref{RRTest}}

The lines bounding the regions defined by \eqref{N1} and \eqref{N3} intersect at $(\frac{1}{p},\frac{1}{s})=(\frac{1}{1+\alpha},\frac{1-\alpha}{1+\alpha})$. Then equality in \eqref{N2} gives $\frac{1}{q}=\frac{1}{1+\alpha}$, so the important estimate is an
$L^{1+\alpha}\rightarrow L^{1+\alpha}(L^{(1+\alpha)/(1-\alpha )})$ estimate. 
Easy estimates combined with an interpolation argument show that Theorem \ref{RRTest} will follow if we establish 
\eqref{est} for $f=\chi_E$ and a collection of triples $(p,q,r)$ which are arbitrarily close to $\big(1+\alpha , 
1+\alpha ,(1+\alpha )/(1-\alpha )\big)$. Standard arguments then show that it is enough to prove that if
$R\chi _E (t,\omega )\geq \mu$ for
$$
(t,\omega )\in F=\{(t,\omega ):\omega\in A ,\, t\in B(\omega )\subset [-1,1]\},
$$
where there is some $B$ such that $B\leq m_1 (B(\omega ))\leq 2B$ for $\omega\in A$, then 
\begin{equation*}
\mu^{p}\lambda (A)^{p/q}\,B^{p/s}\leq C(\delta )\, m_2 (E)
\end{equation*}
if
\begin{equation*}
p=\frac{\alpha +\delta\alpha +1}{\delta\alpha +1},\, q=\alpha +\delta\alpha +1,\, s=\frac{\alpha +\delta\alpha +1}{\delta\alpha +1-\alpha}
\end{equation*}
for small $\delta >0$.

For each $\omega\in A$ let 
\begin{equation*}
E(\omega )=\{t\,\omega +s\,\omega^\perp \in E:t\in B(\omega ),s\in [-1,1] \}.
\end{equation*}
Since $R\chi _E (t,\omega )\geq \mu$ and $m_1 (B(\omega ))\geq B$, it follows that 
\begin{equation}\label{ineq1}
m_2 (E(\omega ))\geq \mu\,B.
\end{equation}
Using the change of coordinates $x\mapsto (x\cdot\omega_1 ,x\cdot\omega_2 )$, one can check that 
\begin{equation}\label{ineq2}
m_2 \big(E(\omega _1 )\cap E(\omega _2 )\big)\lesssim \frac{B^2}{|\omega_1 -\omega _2 |}.
\end{equation}
We will bound $m_2 (E)$ from below by using 
\begin{equation}\label{ineq3}
m_2 (E)\geq m_2 \big(\cup_{j=1}^N E(\omega_j )\big)\geq \sum_{j=1}^N m_2  (E(\omega _j ))-
\sum_{1\leq j<k\leq N}m_2 \big(E(\omega _j )\cap E(\omega _k )\big)
\end{equation}
for appropriately chosen $\omega_j \in A$.
Fix, for the moment, a small positive number $\eta$ 
and consider a partitioning of $S^1$ into intervals of length about $\eta$. Since $\lambda (B(x,r))\lesssim r^\alpha$, 
the $\lambda$-measure of each of these intervals is $\lesssim \eta^{\alpha}$. So at
least, roughly, $\eta^{-\alpha}\lambda (A)$ of them must intersect $A$. Thus it is possible to choose $N\sim \eta^{-\alpha}\lambda (A)$ points $\omega_j \in A$ with $|\omega_j -\omega_k |\gtrsim \eta\, |j-k|$.
Then, for any $\delta >0$, 
$$
\sum_{1\leq j<k\leq N}\frac{1}{|\omega_j -\omega_k |}\lesssim
\eta^{-1} \sum_{1\leq j<k\leq N}\frac{1}{|j -k |}\lesssim \eta^{-1} N^{1+\delta}
$$
and so, by \eqref{ineq2},  
\begin{equation}\label{ineq4}
\sum_{1\leq j<k\leq N}m_2 \big(E({\omega_1})\cap E({\omega_2})\big)
%\lesssim 2^{-N\beta}2^{m} K^{1+\delta}
\leq C\,  B^2 \eta^{-1} N^{1+\delta}\leq C_1 B^2 N^{1+\delta +1/\alpha}\lambda (A)^{-1/\alpha},
\end{equation}
where we have used $N\sim \eta^{-\alpha}\lambda (A)$.
We would now like to choose $N$ such that 
\begin{equation}\label{ineq5}
2\,C_1 \, B^2 N^{1+\delta+1/\alpha}\lambda (A)^{-1/\alpha}\leq N\,\mu \, B\leq
3\,C_1 \, B^2 N^{1+\delta+1/\alpha}\lambda (A)^{-1/\alpha}
\end{equation}
or
\begin{equation}\label{ineq6}
3^{-\alpha/(1+\delta\alpha)}\Big(\frac{\mu\, B^{-1}\lambda (A)^{1/\alpha}}{C_1}\Big)^{\alpha /(\delta\alpha +1)}
\leq N\leq
2^{-\alpha/(1+\delta\alpha)}\Big(\frac{\mu \, B^{-1}\lambda (A)^{1/\alpha}}{C_1}\Big)^{\alpha /(\delta\alpha +1)}.
\end{equation}
This will be possible unless 
$$
\mu \, B^{-1}\lambda (A)^{1/\alpha}\lesssim 1
$$
in which case
$$
\mu^{\alpha /(\delta\alpha +1)}B^{-\alpha /(\delta\alpha +1)}\lambda (A)^{1/(\delta\alpha +1)}\lesssim1
$$
so that the desired inequality  
\begin{equation}\label{ineq7}
m_2 (E)\gtrsim \mu^{(\alpha +\delta\alpha +1)/(\delta\alpha +1)}\lambda (A)^{1/(\delta\alpha +1)}B^{(\delta\alpha +1-\alpha)/(\delta\alpha +1)}
\end{equation}
follows from $m_2 (E)\geq \mu \, B$ unless $F$ is empty. 
Now (with $N$  chosen so that \eqref{ineq5} and \eqref{ineq6} are valid),  \eqref{ineq3}, \eqref{ineq1}, \eqref{ineq4}, and the left member of \eqref{ineq5} give 
$
m_2 (E)\gtrsim N\,\mu\,B.
$
Then the left member of \eqref{ineq6} gives 
\eqref{ineq7} again.

\section{Proof of Theorem \ref{result2}}

For $\rho >0$, let $K_\rho$ be the kernel defined on $\bbR ^d$ by $K_\rho (x)=|x|^{-\rho}\chi_{B(0,R)}(x)$ where 
$R=R(d)$ is positive.
Suppose that the finite nonnegative Borel measure $\nu$ is a $\gamma$-dimensional measure on $\bbR ^d$  
in the sense that $\nu \big( B(x,\delta )\big)\leq C(\nu )\, \delta ^{\gamma}$ for all $x\in\bbR ^d$ 
and $\delta >0$. If $\rho <\gamma$ it follows that 
\begin{equation*}
\nu\ast K_\rho \in L^{\infty}(\bbR ^d ).
\end{equation*}
Also 
\begin{equation*}
\nu\ast K_\rho \in L^{1}(\bbR ^d )
\end{equation*}
so long as $\rho <d$. Thus, for $\epsilon>0$, 
\begin{equation}\label{ineq8}
\nu\ast K_\rho \in L^{p}(\bbR ^d ), \ \rho=\gamma +\frac{1}{p}(d-\gamma)-\epsilon
\end{equation}
by interpolation. 
The following lemma is a weak converse of this observation.
\begin{lemma}\label{besov} If \eqref{ineq8} holds with $\epsilon =0$ and $p>1$, then 
$\nu$ is absolutely continuous with respect to Hausdorff measure of dimension
$\gamma -\epsilon$ for any $\epsilon >0$. Thus the support of $\nu$ has Hausdorff dimension 
at least $\gamma$.
\end{lemma} 

\begin{proof}

Recall from \cite{BL} (see p. 140) that, for $s\in\bbR$ and $1\leq p,q\leq\infty$, the norm $\|f\|^s_{p,q}$ of a distribution $f$ on $\bbR ^d$ in the Besov space $B^s_{p,q}$ can be defined by

$$
\|f\|^s_{pq}=\|\psi\ast f\|_{L^{p}(\bbR ^d )}+\Big(\sum_{k=1}^{\infty} \big(2^{sk}\,\|\phi _k \ast f\|_{L^{p}(\bbR ^d )}\big)^q \Big)^{1/q}
$$
for certain fixed  $\psi\in\cS (\bbR^d )$, $\phi \in C^{\infty}_c (\bbR^d )$, and where $\phi _k (x)=2^{kd}\phi (2^k x)$. If $\nu \ast K_\rho \in L^{p}(\bbR ^d )$, then $\|\nu \ast \chi_{B(0,\delta )} \|_{ L^{p}(\bbR ^d )}\lesssim \delta ^\rho$. It follows that 
$\|\nu\|^s_{pq}<\infty$ if $s<\rho -d=(\gamma -d)/p'$. Now, for $t>0$ and $1<p',q'<\infty$, the Besov capacity $A_{t,p',q'}(K)$ of a compact $K\subset \bbR ^d$ is defined by

$$
A_{t,p',q'}(K)=\inf \{\|f\|^t_{p',q'}:f\in C^{\infty}_c (\bbR^d ),\,f\geq \chi_K \}.
$$
It is shown in \cite{S} (see p. 277) that $A_{t,p',q'}(K)\lesssim H_{d-tp'}(K)$. Thus it follows from the duality of $B^s_{p,q}$ and $B^{-s}_{p',q'}$ that 

$$
\nu (K)\lesssim \|\nu\|^s_{pq}\, A_{-s,p',q'}(K)\lesssim   H_{d+sp'}(K)=H_{\gamma -\epsilon} (K)
$$
if $s=(\gamma -d-\epsilon)/p'$.

\end{proof}

Now suppose that $\mu$ is a nonnegative and compactly supported Borel measure on $\bbR ^2$ which is $\beta$-dimensional in the sense that $\mu \big(B(x,\delta )\big)\lesssim \delta ^\beta$. If the radii $R(1)$ and $R(2)$ (in the definition of $K_\rho$) are chosen so that $R(1)=1$ and $R(2)$ is large enough, depending on the support of $\mu$, then one can verify directly that 
\begin{equation*} 
\mu_\omega \ast K_{(\rho -1)}(t)\lesssim \int_{-2R(2)}^{2R(2)} \mu\ast K_{\rho} \,(t\omega +s\omega^\perp )\,ds .
\end{equation*}
If $p,q,s$ are such that \eqref{est} holds and if $\rho =\beta +(2-\beta )/p-\epsilon$, so that \eqref{ineq8} implies that 
$\mu\ast K_\rho \in L^p (\bbR ^2 )$, then a rescaling of \eqref{est} gives
\begin{equation}\label{ineq10}
\int_{S_1}\|\mu_\omega \ast K_{(\rho -1)}\|^q_{L^s (\bbR )}\,
d\lambda (\omega )<\infty .
\end{equation}
If we could take $(p,q,s)=\big(1+\alpha ,1+\alpha ,(1+\alpha )/(1-\alpha )\big)$ and $\epsilon =0$ then 
\eqref{ineq10} would yield 
\begin{equation*}
\int_{S_1}\|\mu_\omega \ast K_{\tau}\|^{1+\alpha}_{L^{(1+\alpha )/(1-\alpha )}(\bbR )}\,d\lambda (\omega )<\infty
\end{equation*}
with $\tau =(1-\alpha+\alpha\beta )/(1+\alpha )$. Adjusting for the fact that \eqref{ineq10} actually holds only for $(p,q,s)$ 
close to $\big(1+\alpha ,1+\alpha ,(1+\alpha )/(1-\alpha )\big)$ and with $\epsilon >0$, it still follows that 
\begin{equation*}
\int_{S_1}\|\mu_\omega \ast K_{\tau}\|^{1+\alpha -\epsilon}_{L^{(1+\alpha -\epsilon)/(1-\alpha )}(\bbR )}\,d\lambda (\omega )<\infty
\end{equation*}
with  $\tau =(1-\alpha+\alpha\beta )/(1+\alpha )-\epsilon$ for any $\epsilon>0$. 
With $\nu =\mu_\omega$, $p=(1+\alpha -\epsilon)/(1-\alpha )$, and $d=1$, 
Lemma \ref{besov} then shows that, for any $\epsilon >0$, the Hausdorff dimension of $\mu_\omega$'s support exceeds $\beta /2 -\epsilon$ for $\lambda$-almost all $\omega$'s. Since this is true for any $\alpha$-dimensional measure $\lambda$
and for any $\alpha\in (0,1 )$, it follows that $\dim\{\omega\in S^1 :\dim p_\omega (E) <\beta /{2}\}=0 $ as desired.

\end{document}